\documentclass[letterpaper, 10 pt, conference]{ieeeconf}  
\usepackage{multirow}
\usepackage{graphicx}
\usepackage{amsmath}
\usepackage{amsthm}
\usepackage{mathtools}
\usepackage[psamsfonts]{amssymb}
\usepackage{amsxtra}
\usepackage{hyperref}
\usepackage{paralist}
\usepackage{threeparttable}
\usepackage{cleveref}
\usepackage{cite}
\usepackage{color}
\usepackage[justification=centering]{caption}

\renewcommand{\thesection}{\arabic{section}}

\newtheorem{Theorem}{Theorem}[section]
\newtheorem{Lemma}[Theorem]{Lemma}

\theoremstyle{definition}
\newtheorem{assumption}{Assumption}

\IEEEoverridecommandlockouts                              
\title{\LARGE \bf
Convergence rates for distributed stochastic optimization over random networks
}

\author{Dusan~Jakovetic, Dragana~Bajovic, Anit~Kumar~Sahu and Soummya~Kar
\thanks{The work of DJ and DB was supported in part by the EU Horizon 2020 project I-BiDaaS, project number~{780787}. The work of D. Jakovetic was also supported in part by the Serbian Ministry of Education, Science, and Technological Development, grant 174030. The work of AKS and SK was supported in part by National Science Foundation under grant CCF-1513936.}
\thanks{D. Bajovic is with the Faculty of Technical Sciences, University of Novi Sad 21000 Novi Sad, Serbia
        {\tt\small dbajovic@uns.ac.rs}}%
\thanks{D. Jakovetic is with the Department of Mathematics and Informatics, Faculty of Sciences, University of Novi Sad 21000 Novi Sad, Serbia
	{\tt\small djakovet@uns.ac.rs}}
\thanks{A. K. Sahu and S. Kar are with the Department of Electrical and Computer Engineering, Carnegie Mellon University, Pittsburgh, PA 15213
        {\tt\small \{anits,soummyak\}@andrew.cmu.edu}}%
}

\begin{document}

\maketitle
\thispagestyle{empty}
\pagestyle{empty}

\begin{abstract}
We establish the $O(\frac{1}{k})$ convergence rate for distributed stochastic gradient methods that operate over strongly convex costs and random networks. The considered class of methods is standard – each node performs a weighted average of its own and its neighbors’ solution estimates (consensus), and takes a negative step with respect to a noisy version of its local function’s gradient (innovation). The underlying communication network is modeled through a sequence of temporally independent identically distributed (i.i.d.) Laplacian matrices connected on average, while the local gradient noises are also i.i.d. in time, have finite second moment, and possibly unbounded support. We show that, after a careful setting of the consensus and innovations potentials (weights), the distributed stochastic gradient method achieves a (order-optimal) $O(\frac{1}{k})$ convergence rate in the mean square distance from the solution. This is the first order-optimal convergence rate result on distributed strongly convex stochastic optimization when the network is random and/or the gradient noises have unbounded support. Simulation examples confirm the theoretical findings.

\end{abstract}


\section{Introduction}
\label{section-intro}
Distributed optimization and learning algorithms attract
a great interest in recent years,
thanks to their widespread applications
including distributed estimation in networked systems, e.g., \cite{SoummyaGossipEstimation},
distributed control, e.g.,~\cite{BulloBook},
and big data analytics, e.g.,~\cite{daneshmand2015hybrid}.

In this paper, we study distributed stochastic optimization algorithms that operate over random networks and minimize smooth strongly convex costs.
 We consider standard distributed stochastic gradient methods where at each time step,
 each node makes a weighted average of its own and its neighbors' solution estimates,
 and performs a step in the negative direction of its noisy local gradient. The underlying network is
 allowed to be \emph{randomly varying}, similarly to, e.g., the models in~\cite{AsuRandom2,ASU_Math_Prog,randomNesterov}. More specifically,
the network is modeled through a sequence of independent identically distributed (i.i.d.)
 graph Laplacian matrices, where the network is assumed to be connected on average.
 (This translates into the requirement that the algebraic connectivity of the mean Laplacian matrix is strictly positive.) Random network models are highly relevant in, e.g.,
 internet of things (IoT) and cyber physical systems (CPS) applications, like,
 e.g., predictive maintenance and monitoring in industrial manufacturing
 systems, monitoring smart buildings, etc. Therein,
 networked nodes often communicate through unreliable/intermittent wireless links,
 due to, e.g., low-power transmissions or harsh environments.

The main contributions of the paper are as follows.
We show that, by carefully designing the
consensus and the gradient weights (potentials),
the considered distributed stochastic gradient algorithm achieves the order-optimal $O(1/k)$ rate of decay of the mean squared
distance from the solution (mean squared error -- MSE).
 This is achieved for  twice continuously differentiable strongly convex local costs,
 assuming also that the noisy gradients are unbiased estimates of the true gradients and
 that the noise in gradients has bounded second moment.
  To the best of our knowledge, this is the first time an order-optimal
 convergence rate for distributed strongly convex
 stochastic optimization has been established for random networks.

We now briefly review the literature to help us
contrast this paper from prior work.
 In the context of the extensive literature on distributed
 optimization, the most relevant to our work
 are the references on: 1)
 distributed strongly convex stochastic (sub)gradient methods;
 and 2) distributed (sub)gradient methods over random networks (both
 deterministic and stochastic methods).
 For the former thread of works,
 several papers give explicit
 convergence rates under different assumptions.
 Regarding the underlying network,
 references~\cite{RabbatDistributedStronglyCVX,SayedStochasticOpt}
 consider static networks, while the works~\cite{DistributedMirrorDescent,Kozat,NedicStochasticPush}
  consider deterministic time-varying networks.

 References~\cite{RabbatDistributedStronglyCVX,SayedStochasticOpt}
 consider distributed strongly
 convex optimization for static networks, assuming that the
 data distributions that
 underlie each node's local cost function are equal
 (reference \cite{RabbatDistributedStronglyCVX}
 considers empirical risks while
 reference \cite{SayedStochasticOpt} considers risk functions
 in the form of expectation);
 this essentially corresponds to each nodes'
 local function having the same minimizer.
 References \cite{DistributedMirrorDescent,Kozat,NedicStochasticPush}
   consider deterministically varying networks,
   assuming that the ``union graph'' over
   finite windows of iterations is connected.
    The papers~\cite{RabbatDistributedStronglyCVX,SayedStochasticOpt,DistributedMirrorDescent,Kozat}
     assume undirected networks,
     while \cite{NedicStochasticPush} allows for
     directed networks but assumes a bounded support for the gradient noise.
      The works~\cite{RabbatDistributedStronglyCVX,DistributedMirrorDescent,Kozat,NedicStochasticPush}
       allow the local costs to be non-smooth,
       while \cite{SayedStochasticOpt} assumes smooth costs,
       as we do here.
 With respect to these works, we consider random networks, undirected networks,
  smooth costs, and allow the noise to have unbounded support.

Distributed optimization over random networks has been studied in
\cite{AsuRandom2,ASU_Math_Prog,randomNesterov}. References \cite{AsuRandom2,ASU_Math_Prog}
consider non-differentiable convex costs and no (sub)gradient noise,
 while reference \cite{randomNesterov}
 considers differentiable costs with Lipschitz continuous and bounded
 gradients, and it also does not allow for gradient noise, i.e.,
 it considers methods with exact (deterministic) gradients.
 In~\cite{CDCKW2018},
 we consider a distributed Kiefer-Wolfowitz-type
 stochastic approximation method and establish
 the method's $O(1/k^{1/2})$ convergence rate.
  Reference~\cite{CDCKW2018}
 complements the current paper by assuming
 that nodes only have access to noisy
 functions' estimates (zeroth
 order optimization), and no gradient estimates are available.
 In contrast, by assuming a noisy first-order (gradient) information is
 available, we show here that strictly faster rates than
 in~\cite{CDCKW2018} can be achieved.

In summary, to the best of our knowledge,
this is the first time to establish an order-optimal
convergence rate for distributed stochastic gradient methods
under strongly convex costs and random networks.


\textbf{Paper organization}. The next
paragraph introduces notation.
Section~2 describes the model and the stochastic gradient method we consider.
Section~3 states and proves the main result on the algorithm's MSE convergence rate.
Section~4 provides a simulation example. Finally, we conclude in Section~5.

\textbf{Notation}. We denote by $\mathbb R$ the set of real numbers and by ${\mathbb R}^m$ the $m$-dimensional
Euclidean real coordinate space. We use normal lower-case letters for scalars,
lower case boldface letters for vectors, and upper case boldface letters for
matrices. Further, we denote by: $\mathbf{A}_{ij}$ the entry in the $i$-th row and $j$-th column of
a matrix $\mathbf{A}$;
$\mathbf{A}^\top$ the transpose of a matrix $\mathbf{A}$; $\otimes$ the Kronecker product of matrices;
$\mathbf{I}$, $0$, and $\mathbf{1}$, respectively, the identity matrix, the zero matrix, and the column vector with unit entries; $\mathbf{J}$ the $N \times N$ matrix $J:=(1/N)\mathbf{1}\mathbf{1}^\top$.
When necessary, we indicate the matrix or vector dimension through a subscript.
 Next, $A \succ  0 \,(A \succeq  0 )$ means that
the symmetric matrix $A$ is positive definite (respectively, positive semi-definite).
We further denote by:
 $\|\cdot\|=\|\cdot\|_2$ the Euclidean (respectively, spectral) norm of its vector (respectively, matrix) argument; $\lambda_i(\cdot)$ the $i$-th smallest eigenvalue; $\nabla h(w)$ and $\nabla^2 h(w)$ the gradient and Hessian, respectively, evaluated at $w$ of a function $h: {\mathbb R}^m \rightarrow {\mathbb R}$, $m \geq 1$; $\mathbb P(\mathcal A)$ and $\mathbb E[u]$ the probability of
an event $\mathcal A$ and expectation of a random variable $u$, respectively.
Finally, for two positive sequences $\eta_n$ and $\chi_n$, we have: $\eta_n = O(\chi_n)$ if
 $\limsup_{n \rightarrow \infty}\frac{\eta_n}{\chi_n}<\infty$.

\section{Model and algorithm}
\subsection{Optimization and network models}
We consider the scenario where $N$
networked nodes aim to collaboratively solve the following unconstrained problem:
\begin{align}
\label{eq:opt_problem}
\mathrm{minimize\,\,\,\,}\sum_{i=1}^{N}f_{i}(\mathbf{x}),
\end{align}
where $f_{i}: \mathbb{R}^{m}\mapsto\mathbb{R}$ is a convex
function available to node $i$, $i=1,...,N$.
We make the following standard assumption on the $f_i$'s.
\begin{assumption}
\label{assumption-f-i}
For all $i=1,...,N$, function
$f_{i}: \mathbb{R}^{m}\mapsto\mathbb{R}$ is twice
continuously differentiable, and there exist
constants $0<\mu \leq L < \infty$, such that,
for all $\mathbf{x} \in {\mathbb R}^m$, there holds:
\[
\mu\,\mathbf{I} \preceq \nabla^2 f_i(\mathbf{x}) \preceq L\,\mathbf{I}.
\]
\end{assumption}
Assumption~\ref{assumption-f-i} implies that
each $f_i$ is strongly convex with
modulus $\mu$, and it also
has Lipschitz continuous
gradient with Lipschitz constant~$L$,
i.e., the following two inequalities hold
 for any $\mathbf{x,y} \in {\mathbb R}^m$:
\begin{eqnarray*}
&\,&f_i(\mathbf{y}) \geq
f_i(\mathbf{x})
+ \nabla f_i(\mathbf{x})^\top \,(\mathbf{y}-\mathbf{x})
+ \frac{\mu}{2}\|\mathbf{x}-\mathbf{y}\|^2\\
&\,&\|\nabla f_i(\mathbf{x}) -
\nabla f_i(\mathbf{y})\| \leq
L\,\|\mathbf{x}-\mathbf{y}\|.
\end{eqnarray*}
Furthermore,
under Assumption~\ref{assumption-f-i},
 problem \eqref{eq:opt_problem} is solvable
 and has the unique solution, which we
 denote by $x^\star \in {\mathbb R}^m$.
 For future reference, introduce also the sum function
 $f:\,{\mathbb R}^m \rightarrow \mathbb R$,
 $f(\mathbf{x})
 =\sum_{i=1}^N f_i(\mathbf{x})$.

We consider distributed stochastic
gradient methods to solve~\eqref{eq:opt_problem}
over random networks. Specifically,
we adopt the following model.
At each time instant $k=0,1,...$,
 the underlying network
 $\mathcal{G}(k) = (V,\mathbf{E}(k))$
 is undirected and random,
 with $V=\{1,...,N\}$
 the set of nodes,
 and $\mathbf{E}(k)$
 the random set of undirected
 edges. We denote by $\{i,j\}$
 the edge that connects nodes $i$ and~$j$.
 Further, denote by 7
 $\Omega_i(k)=\{j \in V:\,\,\{i,j\} \in \mathbf{E}(k)\}$
  the random neighborhood od
  node $i$ at time $k$ (excluding node~$i$).
 We associate to $\mathcal{G}(k)$ its
 $N \times N$ (symmetric) Laplacian matrix
 $\boldsymbol{\mathcal{L}}(k)$, defined
 by:
 $
  \boldsymbol{\mathcal{L}}_{ij}(k) =
  -1,$ if $\{i,j\} \in \mathbf{E}(k)$, $i \neq j$;
  $
  \boldsymbol{\mathcal{L}}_{ij}(k) =
  0,$ if $\{i,j\} \notin \mathbf{E}(k)$, $i \neq j$;
  and~$\boldsymbol{\mathcal{L}}_{ii}(k) =
  -\sum_{j \neq i}\boldsymbol{\mathcal{L}}_{ij}(k).$
  Denote by $\overline{\boldsymbol{\mathcal{L}}}=
  \mathbb{E}[\, \boldsymbol{\mathcal{L}}(k)\,]$ (as explained ahead,
  the expectation is independent of $k$),
  and let $\overline{\mathcal{G}} = (V,\overline{\mathbf{E}})$
   be the graph induced by matrix $\overline{\boldsymbol{\mathcal{L}}}$,
   i.e.,
   $\overline{\mathbf{E}} = \{\{i,j\}:\,\,i \neq j,\,\,\,
   \overline{\boldsymbol{\mathcal{L}}}_{ij}>0
   \}$.
 We make the following assumption.
\begin{assumption}
\label{assumption-network}
The matrices $\{\mathbf{L}(k)\}$
are independent, identically distributed (i.i.d.)
Furthermore, graph $\overline{\mathcal{G}}$
 is connected.
\end{assumption}
It is well-known that
the connectedness of $\overline{\mathcal{G}}$
is equivalent to
the condition $\lambda_2
\left(
\overline{\boldsymbol{\mathcal{L}}} \right)>0$.

\subsection{Gradient noise model and the algorithm}
We consider the following
distributed stochastic
gradient method to solve~\eqref{eq:opt_problem}.
Each node~$i$, $i=1,...,N$,
maintains over time steps
(iterations) $k=0,1,...,$
 its solution estimate
 $\mathbf{x}_i(k) \in {\mathbb R}^m$.
 Specifically, for
 arbitrary deterministic initial points
 $\mathbf{x}_i(0) \in {\mathbb R}^m$, $i=1,...,N$,
 the update rule at node $i$
 and $k=0,1,...,$
 is as follows:
 \begin{eqnarray}
 \label{eqn-alg-node-i}
 \mathbf{x}_i(k+1) &=&
 \mathbf{x}_i(k) - \beta_k \,
 \sum_{j \in \Omega_i(k)}
  \left( \mathbf{x}_i(k) - \mathbf{x}_j(k) \right)\\
  &-&
  \alpha_k\,\left( \,\nabla f_i(\mathbf{x}_i(k)) + \mathbf{v}_i(k)\,\right).\nonumber
 \end{eqnarray}
The update \eqref{eqn-alg-node-i} is
performed in parallel by all nodes $i=1,...,N$.
The algorithm iteration is realized as follows. First,
each node $i$
broadcasts $\mathbf{x}_i(k)$ to
all its available neighbors
$j \in \Omega_i(k)$,
and receives
$\mathbf{x}_j(k)$ from all
$j \in \Omega_i(k)$.
Subsequently,
each node $i$, $i=1,...,N$
 makes update~\eqref{eqn-alg-node-i},
 which completes an iteration.
 In \eqref{eqn-alg-node-i},
 $\alpha_k$
  is the step-size
  that we set to
  $\alpha_k=\alpha_0/(k+1)$,
  $k=0,1,...,$ with $\alpha_0>0$;
  and $\beta_k$ is
  the (possibly) time-varying weight
  that each node
  assigns to all its neighbors.
  We set
   $\beta_k=\beta_0/(k+1)^{\nu}$,
  $k=0,1,...$, with
  $\nu \in [0,1/2]$. Here,
  $\beta_0>0$ is a constant
  that should be taken to be sufficiently small; e.g.,
  one can set $\beta_0=1/(1+
  \theta)$, where $\theta$ is the maximal degree (number of neighbors of a node)
  across network.
  Finally,
  $\mathbf{v}_i(k)$
  is noise
  in the calculation
  of the $f_i$'s gradient at iteration~$k$.

For future reference, we also
present algorithm \eqref{eqn-alg-node-i} in matrix format.
 Denote by $\mathbf{x}(k) = \left[\mathbf{x}_{1}^{\top}(k),\cdots,\mathbf{x}_{N}^{\top}(k)\right]^{\top}\in\mathbb{R}^{Nm}$
  the vector that stacks the solution estimates of all nodes.
  Also, define function $F: \mathbb{R}^{Nm}\mapsto\mathbb{R}$, by
  $F(\mathbf{x})=\sum_{i=1}^{N}f_{i}(\mathbf{x}_{i})$,
   with $\mathbf{x} = \left[\mathbf{x}_{1}^{\top},\cdots,\mathbf{x}_{N}^{\top}\right]^{\top}\in\mathbb{R}^{Nm}.$
    Finally, let $\mathbf{W}_{k}=\left(\mathbf{I}-\mathbf{L}_{k}\right)\otimes \mathbf{I}_{m}$, where $\mathbf{L}_{k}=\beta_k\,\boldsymbol{\mathcal{L}}(k)$.
    Then, for $k=0,1,...$, algorithm
    \eqref{eqn-alg-node-i} can be compactly written as follows:
\begin{align}
\label{eq:update_rule}
\mathbf{x}(k+1) = \mathbf{W}_{k}\mathbf{x}(k)-\alpha_{k}\left(\nabla F(\mathbf{x}(k))+\mathbf{v}(k)\right).
\end{align}


  We make the following standard
  assumption on the gradient noises.
  First, denote by $\mathcal{F}_k$
   the history
   of algorithm \eqref{eqn-alg-node-i}
   up to time $k$;
   that is, $\mathcal{F}_k$, $k=1,2,...,$
   is an increasing sequence of sigma algebras,
    where
    $\mathcal{F}_k$
     is the sigma algebra generated
     by the collection of random variables $
     \{ \,\boldsymbol{\mathcal{L}}(s),
     \,\mathbf{v}_i(t)\}$,
     $i=1,...,N$,
     $s=0,...,k-1$,
     $t=0,...,k-1$.
   \begin{assumption}
   \label{assumption-gradient-noise}
   For each $i=1,...,N$,
   the sequence of noises $\{\mathbf{v}_i(k)\}$
   satisfies for all $k=0,1,...$:
   \begin{eqnarray}
   \label{eqn-ass-noise-1}
    \mathbb{E}[\,\mathbf{v}_i(k)\,|\,\mathcal{F}_k\,] &=& 0,\,\,\mathrm{almost\,surely\,(a.s.)}\\
   \label{eqn-ass-noise-2}
    \mathbb{E}[\,\|\mathbf{v}_i(k)\|^2\,|\,\mathcal{F}_k\,] &\leq& c_{v}\|\mathbf{x}_i(k)\|^2+c_v^\prime,
    \,\,\mathrm{a.s.},
   \end{eqnarray}
   where $c_v$ and $c_v^\prime$
   are nonnegative constants.
   \end{assumption}
   %
   %
   %
%
Assumption~\ref{assumption-gradient-noise} is satisfied,
for example, when $\{\mathbf{v}_i(k)\}$
is an i.i.d. zero-mean, finite second moment, noise sequence such that
 $\mathbf{v}_i(k)$ is also independent
 of history~$\mathcal{F}_k$.
     However, the assumption
     allows that
     the gradient noise
     $\mathbf{v}_i(k)$
     be dependent on node $i$
      and also on the
      current point
      $\mathbf{x}_i(k)$; the next subsection
      gives some important machine learning settings encompassed by
      Assumption~\ref{assumption-gradient-noise}.

\subsection{A machine learning motivation}
The optimization-algorithmic model
defined by Assumptions
\ref{assumption-f-i} and \ref{assumption-gradient-noise}
subsumes, e.g., important machine learning applications.
 Consider the scenario
 where $f_i$ corresponds to
 the risk function associated
 with the node $i$'s local data, i.e.,
 \begin{equation}
 \label{eqn-risk-fcn}
 f_i(\mathbf{x}) = \mathbb{E}_{\,\mathbf{d}_i \sim P_i}\left[\,\ell_i\left(\mathbf{x};\mathbf{d}_i\right) \,\right]+\Psi_i(\mathbf{x}).
 \end{equation}
Here, $P_i$ is node $i$'s local distribution according to
which its data samples $\mathbf{d}_i \in {\mathbb R}^q$ are generated; $\ell_i(\cdot;\cdot)$
is a loss function
that is convex in
its first argument
for any fixed value of its second argument;
and $\Psi:\,{\mathbb R}^m \rightarrow \mathbb R$
 is a strongly convex regularizer.
 Similarly,
 $f_i$ can be an empirical risk function:
 \begin{equation}
 \label{eqn-empir-risk-fcn}
 f_i(\mathbf{x}) = \frac{1}{n_i}\left(\,\sum_{j=1}^{n_i}
 \ell_i\left(\mathbf{x};\mathbf{d}_{i,j}\right)\,\right) +\Psi_i(\mathbf{x}),
 \end{equation}
 where $\mathbf{d}_{i,j}$,
 $j=1,...,n_i$, is the set of training examples at node~$i$.
  Examples for
  the loss $\ell_i(\cdot;\cdot)$ include the following:
\begin{eqnarray}
\label{eqn-ML-examples}
\ell_i(\mathbf{x};\mathbf{a_i},b_i)
&=&\frac{1}{2} \left( \mathbf{a}_i^\top \mathbf{x} - b_i\right)^2 \,\,\,\,\,(\mathrm{quadratic\,\,loss})\\
\ell_i(\mathbf{x};\mathbf{a_i},b_i)
&=&
\mathrm{ln}\left( \,1+\mathrm{exp}(-b_i(\mathbf{a}_i^\top \mathbf{x}))\,\right)
\,\,\,\,\,(\mathrm{logistic\,\,loss}) \nonumber
\end{eqnarray}
For the quadratic loss above,
a data sample $\mathbf{d}_i=(\mathbf{a}_i,b_i)$,
where $\mathbf{a}_i$ is a regressor vector and
$b_i$ is a response variable;
for the logistic loss,
 $\mathbf{a}_i$ is a feature vector and
 $b_i \in \{-1,+1\}$ is its class label.
 Clearly,
 both the risk \eqref{eqn-risk-fcn}
 and the empirical risk \eqref{eqn-empir-risk-fcn}
 satisfy Assumption~\ref{assumption-f-i} for
 the losses in~\eqref{eqn-ML-examples}.

  We next discuss the search directions
  in \eqref{eqn-alg-node-i} and Assumption~\ref{assumption-gradient-noise} for the
  gradient noise.
  A common
  search direction in machine learning
  algorithms is
  the gradient of the loss with respect to
  a single data point\footnote{Similar considerations hold for
  a loss with respect to a mini-batch of data points; this discussion is
  abstracted for simplicity.}:
  \[
  g_i(\mathbf{x}) = \nabla \ell_i\left(\mathbf{x};\mathbf{d}_i\right) + \nabla \Psi_i(\mathbf{x}).
  \]
In case of the risk function \eqref{eqn-ML-examples},
$\mathbf{d}_i$ is drawn from distribution~$P_i$;
in case of the empirical risk \eqref{eqn-empir-risk-fcn},
$\mathbf{d}_i$ can be, e.g., drawn uniformly
at random from the set of
data points $\mathbf{d}_{i,j}$, $j=1,...,n_i$,
with repetition along iterations.
In both cases,
gradient noise
 $\mathbf{v}_i = g_i(\mathbf{x}) - \nabla f_i(\mathbf{x})$
 clearly satisfies
assumption~\eqref{eqn-ass-noise-1}.
 To see this, consider, for example,
the risk function \eqref{eqn-risk-fcn},
and let us fix iteration $k$
 and node $i$'s estimate $\mathbf{x}_i(k) = \mathbf{x}_i$.
  Then,
\begin{eqnarray*}
&\,&\mathbb{E}
\left[ \mathbf{v}_i(k)\,|\,\mathcal{F}_k\right] =
\mathbb{E}
\left[ \mathbf{g}_i(k)-\nabla f_i(\mathbf{x}_i(k))\,|\,\mathbf{x}_i(k) = \mathbf{x}_i\right]\\
&=&
\mathbb{E} [\nabla \ell_i\left(\mathbf{x}_i(k);\mathbf{d}_i\right)\,|\,\mathbf{x}_i(k)=\mathbf{x}_i]
 + \nabla \Psi_i(\mathbf{x}_i)\\
 &-&
 \left( \nabla f_i(\mathbf{x}_i)
 + \nabla \Psi_i(\mathbf{x})\right)\\
 &=&
\mathbb{E}_{\,\mathbf{d}_i \sim P_i}[\nabla \ell_i\left(\mathbf{x};\mathbf{d}_i\right)]
 + \nabla \Psi_i(\mathbf{x}_i)\\
 &-&
 \left(\mathbb{E}_{\,\mathbf{d}_i \sim P_i}[\nabla \ell_i\left(\mathbf{x};\mathbf{d}_i\right)]
 + \nabla \Psi_i(\mathbf{x}_i)\right)=0.
\end{eqnarray*}
Further, for the
empirical risk,
assumption~\eqref{eqn-ass-noise-2}
 holds trivially.
  For the risk function \eqref{eqn-risk-fcn},
assumption~\eqref{eqn-ass-noise-2} holds for
a sufficiently ``regular'' distribution~$P_i$.
 For instance, it is easy to show that the assumption
holds for the logistic loss in~\eqref{eqn-ML-examples}
when $P_i$ has finite second moment,
while it holds for the square loss in~\eqref{eqn-ML-examples}
when $P_i$ has finite fourth moment.

Note that our setting allows
that the data generated at different nodes
be generated through different distributions~$P_i$,
as well as that the nodes utilize different losses
$\ell_i$'s and regularizers $\Psi_i$'s.
 Mathematically, this means
 that $\nabla f_i(x^\star) \neq 0$, in general.
 In words, if a node~$i$
 relies only on its local data $\mathbf{d}_i$,
 it cannot recover the true solution~$x^\star$.
  Nodes then engage in a collaborative
  algorithm~\eqref{eqn-alg-node-i}
   through which, as shown ahead, they can recover
   the global solution~$x^\star$.

\section{Performance Analysis}
\label{sec:perf_analysis}

\subsection{Statement of main results and auxiliary lemmas}
We are now ready to state our main result.
\begin{Theorem}
\label{theorem-1}
Consider algorithm \eqref{eqn-alg-node-i} with
step-sizes $\alpha_k=\frac{\alpha_0}{k+1}$
 and $\beta_k=\frac{\beta_0}{(k+1)^{\nu}}$,
  where $\beta_0>0$,
  $\alpha_0>2\,N/\mu$, and
  $\nu \in [0,1/2]$. Further, let Assumptions 1--3 hold.
Then, for each node~$i$'s
solution estimate $\mathbf{x}_i(k)$ and
the solution $\mathbf{x}^\star$ of problem~\eqref{eq:opt_problem}, there holds:
\[
\mathbb{E}\left[ \|\mathbf{x}_i(k) - \mathbf{x}^\star\|^2\right]=O(1/k).
\]
\end{Theorem}
We remark that
the condition
$\alpha_0>2\,N/\mu$ can be
relaxed to
require only a positive
$\alpha_0$,
in which case the rate becomes
$O(\mathrm{ln}(k)/k)$, instead
of~$O(1/k)$.\footnote{This subtlety
comes from equation~(32) ahead
and the requirement that $c_{20}>1$.
If $c_{20} \leq 1$, it can be shown
that in~\eqref{eqn-invoking}
 the right hand side
 modifies to a $O(\mathrm{ln}(k)/k)$ quantity.}
 Also, to avoid
 large step-sizes
 at initial iterations for a large $\alpha_0$,
 step-size $\alpha_k$
 can be modified to
 $\alpha_k=\alpha_0/(k+k_0)$,
 for arbitrary positive constant~$k_0$,
 and Theorem~\ref{theorem-1}
 continues to hold.
Theorem~\ref{theorem-1}
establishes the $O(1/k)$ MSE rate of convergence
of algorithm~\eqref{eqn-alg-node-i}; due to the
assumed $f_i$'s strong convexity, the theorem also implies
that $\mathbb{E}\left[ f(\mathbf{x}_i(k)) - f(\mathbf{x}^\star)\right]
=O(1/k)$. Note that the expectation
in Theorem~\ref{theorem-1} is both
with respect to randomness in gradient noises
and with respect to the randomness in the underlying network.
 The $O(1/k)$ rate does not depend
 on the statistics of the underlying random network,
 as long as the network is connected on average (i.e.,
 satisfies Assumption~\ref{assumption-network}.)
 The hidden constant depends on
 the underlying network statistics,
 but simulation examples
 suggest that the dependence
 is usually not strong
 (see Section~4).

\textbf{Proof strategy and auxiliary lemmas}.
Our strategy for proving Theorem~\ref{theorem-1}
is as follows. We first
establish the mean square boundedness (uniform in~$k$)
 of the iterates $\mathbf{x}_i(k)$, which also
 implies the uniform mean square
 boundedness of the gradients $\nabla f_i(\mathbf{x}_i(k))$
 (Subsection~3-B).
 We then bound, in the mean square sense, the
 disagreements of
 different nodes' estimates, i.e.,
 quantities~$(\mathbf{x}_i(k)-\mathbf{x}_j(k))$, showing that
 $\mathbb{E}[\,\|\mathbf{x}_i(k)-\mathbf{x}_j(k)\|^2\,]=O(1/k)$
 (Subsection~3-C).
 This allows us to show that
 the (hypothetical)
 global average of
 the nodes'
 solution estimates
 $\overline{\mathbf{x}}(k):=\frac{1}{N}\sum_{i=1}^N \mathbf{x}_i(k)$
 evolves according to
 a stochastic gradient method
 with the gradient estimates
 that have a sufficiently small bias and finite second moment.
  This allows us to show the
  $O(1/k)$ rate on
  the mean square error at the global average, which in turn
  allows to derive a similar bound at the individual nodes' estimates
  (Subsection~3-D).

In completing the strategy above, we make use of the following Lemma;
the Lemma is a minor modification of Lemmas~4 and~5 in~\cite{SoummyaGossipEstimation}.
\begin{Lemma}
\label{lemma-estimation}
Let $z(k)$ be a nonnegative (deterministic) sequence satisfying:
\[
z(k+1) \leq (1-r_1(k))\,z_1(k) + r_2(k),
\]
where $\{r_1(k)\}$ and $\{r_2(k)\}$
are deterministic sequences with
\begin{eqnarray*}
\frac{a_1}{(k+1)^{\delta_1}}
\leq r_1(k) \leq 1 \,\,\,\mathrm{and}\,\,\,
r_2(k) \leq
\frac{a_2}{(k+1)^{\delta_2}},
\end{eqnarray*}
with $a_1 , a_2 , \delta_1 , \delta_2 > 0.$ Then, (a) if $\delta_1 =
\delta_2 = 1$, there holds:
$z(k)=O(1)$; (b) if $\delta_1=1/2$ and $\delta_2=3/2$, then $z(k)=O(1/k)$;
and (c) if $\delta_1=1$, $\delta_2=2$, and $a_1>1$, then~$z(k)=O(1/k)$.
\end{Lemma}
%
%
%
%

Subsequent analysis in Subsections~3-b until 3-d restricts
to the case when $\nu=1/2$, i.e.,
when consensus weights equal
$\beta_k=\frac{\beta_0}{(k+1)^{1/2}}$. That is,
for simplicity of presentation,
we prove Theorem~\ref{theorem-1}
for case $\nu=1/2$.
As it can be verified in subsequent
analysis, the proof of Theorem~\ref{theorem-1} extends
to a generic $\mu \in [0,1/2)$ as well.
 As another step in
 simplifying notations, throughout Subsections~3-b and 3-c,
we let $m=1$ to avoid extensive usage of Kronecker products; again,  the proofs extend to a generic $m>1.$

\subsection{Mean square boundedness of the iterates}
%
This Subsection
shows the uniform mean square boundedness
of the algorithm iterates and the gradients evaluated at the algorithm
iterates.
\begin{Lemma}
\label{Lemma-MSS-BDD}
Consider algorithm \eqref{eqn-alg-node-i}, and
let Assumptions~1-3 hold.
Then, there exist
nonnegative constants
$c_x$ and $c_{\,\partial f}$
such that, for all $k=0,1,...,$ there holds:
\[
\mathbb{E}[\,\|\mathbf{x}(k)\|^2\,] \leq c_x\,\,\,\,
\mathrm{and}\,\,\,\,
\mathbb{E}[\,\|\nabla F(\mathbf{x}(k))\|^2\,] \leq c_{\,\partial f}.
\]
\end{Lemma}

\textit{Proof.}
%

Denote by $\mathbf{x}^{o} = x^{\ast}\mathbf{1}_{N}$ and recall \eqref{eq:update_rule}.
 Then, we have:
\begin{eqnarray}
\label{eq:update_rule2}
\mathbf{x}(k+1)-\mathbf{x}^{o} &=& \mathbf{W}_{k}(\mathbf{x}(k)-\mathbf{x}^{o})\\
&-&\alpha_{k}\left(\nabla F(\mathbf{x}(k))-\nabla F(\mathbf{x}^{o})\right) \nonumber \\
&-&\alpha_{k}\mathbf{v}(k)-\alpha_{k}\nabla F(\mathbf{x}^{o}). \nonumber
\end{eqnarray}
By mean value theorem, we have:
\begin{eqnarray}
\label{eq:mvt}
&\,& \nabla F(\mathbf{x}(k))-\nabla F(\mathbf{x}^{o}) \\
&=& \left[\int_{s=0}^{1}\nabla^{2}F\left(\mathbf{x}^{o}+s({\mathbf{x}(k)}-
\mathbf{x}^{o})\right)\,d \,s\right]\left(\mathbf{x}(k)-\mathbf{x}^{o}\right)\nonumber\\
&=&\mathbf{H}_{k}\left(\mathbf{x}(k)-\mathbf{x}^{o}\right). \nonumber
\end{eqnarray}
Note that $L\mathbf{I}\succcurlyeq\mathbf{H}_{k}\succcurlyeq\mu\mathbf{I}$. Using \eqref{eq:mvt} in  \eqref{eq:update_rule2} we have:
\begin{eqnarray}
\label{eq:update_rule3}
&\,& \mathbf{x}(k+1)-\mathbf{x}^{o} = \left(\mathbf{W}_{k}-\alpha_{k}\mathbf{H}_{k}\right)(\mathbf{x}(k)-\mathbf{x}^{o})\\
&-& \alpha_{k}\mathbf{v}(k)-\alpha_{k}\nabla F(\mathbf{x}^{o}).\nonumber
\end{eqnarray}
%
Denote by
$\boldsymbol{\zeta}(k) = \mathbf{x}(k)-\mathbf{x}^{o}$
and by $\boldsymbol{\xi}(k) = \left(\mathbf{W}_{k}-\alpha_{k}\mathbf{H}_{k}\right)(\mathbf{x}(k)-\mathbf{x}^{o})
-\alpha_k \,\nabla F(\mathbf{x}^{o})$.
Then, there holds:
\begin{eqnarray}
&\,&\mathbb{E}[\,\|\boldsymbol{\zeta}(k+1)\|^2 \,|\,\mathcal{F}_k\,]
\leq
\|\boldsymbol{\xi}(k)\|^2 \nonumber \\
&-& 2 \alpha_k \,{\boldsymbol{\xi}}(k)^\top
\mathbb{E}[\,\mathbf{v}(k) \,|\,\mathcal{F}_k\,] +
\alpha_k^2 \,\mathbb{E}[\,\|\mathbf{v}(k)\|^2 \,|\,\mathcal{F}_k\,] \nonumber \\
&\leq&
\|\boldsymbol{\xi}(k)\|^2 + N\,\alpha_k^2\,(c_v\,\|\mathbf{x}(k)\|^2+c_v^\prime), \,\,\mathrm{a.s.},
\label{eqn-combine-1}
\end{eqnarray}
where we used Assumption~\ref{assumption-gradient-noise}
 and the fact that $\boldsymbol{\xi}(k)$
  is measurable with respect to $\mathcal{F}_k$.
We next bound $\|\boldsymbol{\xi}(k)\|^2$.
 Note that
 $\|\mathbf{W}_k -\alpha_k \,\boldsymbol{H}_k\| \leq 1-\mu\,\alpha_k$ for sufficiently large $k$.
 Therefore, we have for sufficiently large $k$:
 \begin{equation}
 \label{eqn-xi-zeta}
 \|\boldsymbol{\xi}(k)\| \leq (1-\mu\,\alpha_k)\,\|\boldsymbol{\zeta}(k)\|
 + \alpha_k\, \|\nabla F(\mathbf{x}^{o})\|.
 \end{equation}
We now
use the following inequality:
\begin{align}
\label{eq:cool_ineq}
(a+b)^{2} \leq \left(1+\theta\right)a^{2}+\left(1+\frac{1}{\theta}\right)b^{2},
\end{align}
for any $a,b \in \mathbb{R}$ and $\theta > 0$.
We set $\theta=\frac{c_0}{k+1}$, with $c_0>0$.
Using the inequality \eqref{eq:cool_ineq} in \eqref{eqn-xi-zeta}, we have:
\begin{eqnarray*}
&\,&\left\|\boldsymbol{\xi}(k)\right\|^{2}
\le \left(1+\frac{c_0}{k+1}\right)(1-\alpha_k\mu)^{2}\\
&\times&\left\|\boldsymbol{\zeta}(k)\right\|^2+\left(1+\frac{k+1}{c_0}\right)\alpha_{k}^{2}\|\nabla F(\mathbf{x}^{o})\|^2. \nonumber
\end{eqnarray*}
Next, for $c_0 < \alpha_0\mu$,
the last inequality implies:
\begin{align}
\label{eqn-combine-2}
&\left\|\boldsymbol{\xi}(k)\right\|^{2} \le \left(1-\frac{c_1}{k+1}\right)\left\|\boldsymbol{\zeta}(k)\right\|^2\\
&+\frac{c_2}{k+1}\|\nabla F(\mathbf{x}^{o})\|^2,\nonumber
\end{align}
for some constants~$c_1,c_2>0.$
%
%
%
%
%
Combining
\eqref{eqn-combine-2} and \eqref{eqn-combine-1}, we get:
\begin{eqnarray}
\mathbb{E}[\,\|\boldsymbol{\zeta}(k+1)\|^2 \,|\,\mathcal{F}_k\,]
&\leq&
\left( 1-\frac{c_1^\prime}{k+1} \right)\|\boldsymbol{\zeta}(k)\|^2 \nonumber \\
&+&
\frac{c_2^\prime}{k+1},
\label{eqn-proof-100}
\end{eqnarray}
for some $c_1^\prime, c_2^\prime>0.$
Taking expectation in \eqref{eqn-proof-100}
and applying Lemma \ref{lemma-estimation},
it follows that
$\mathbb{E}[\,\|\boldsymbol{\zeta}(k)\|^2\,]
= \mathbb{E}[\,\|\mathbf{x}(k)-\mathbf{x}^{o}\|^2\,]$
is uniformly (in~$k$) bounded from above by a positive constant. It is easy to see that
the latter implies that  $
\mathbb{E}[\,\|\mathbf{x}(k)\|^2\,]$
is also uniformly bounded.
Using the Lipschitz continuity of
$\nabla F$, we finally also have
that $
\mathbb{E}[\,\|\nabla F(\mathbf{x}(k))\|^2\,]$
is also uniformly bounded. The proof of Lemma \ref{Lemma-MSS-BDD}
is now complete.

%
%

\subsection{Disagreement bounds}
\label{sec:dis_bound}
Recall the (hypothetically available) global average of nodes' estimates
$\overline{\mathbf{x}}(k)=\frac{1}{N}\sum_{i=1}^N \mathbf{x}_i(k)$,
and denote by
$\widetilde{\mathbf{x}}_i(k)
=\mathbf{x}_i(k) - \overline{\mathbf{x}}(k)$ the quantity
that measures how far apart is node $i$'s
solution estimate from the global average.
 Introduce also vector
 $\widetilde{\mathbf{x}}(k)=(\,\widetilde{\mathbf{x}}_1(k),...,\widetilde{\mathbf{x}}_N(k)\,)^\top$,
 and note that it can be represented as  $\widetilde{\mathbf{x}}(k)=\left(\mathbf{I}-\mathbf{J}\right)\mathbf{x}(k)$, where we recall $\mathbf{J}=\frac{1}{N}\mathbf{1}\mathbf{1}^{\top}$.
 We have the following Lemma.
\begin{Lemma}
\label{lemma-disag-bound}
Consider algorithm \eqref{eqn-alg-node-i}
under Assumptions~1--3. Then, there holds:
\[
\mathbb{E}[\,\|\widetilde{\mathbf{x}}(k)\|^2\,] = O(1/k).
\]
\end{Lemma}
As detailed in the next Subsection,
Lemma~\ref{lemma-disag-bound}
 is important as it allows
 to sufficiently tightly
 bound the bias in the
 gradient estimates
 according to which the
 global average
 $\overline{\mathbf{x}}(k)$ evolves.

\textit{Proof.}
 It is easy to show that the process $\{\widetilde{\mathbf{x}}(k)\}$ follows the recursion:
\begin{align}
\label{eq:dis1}
\widetilde{\mathbf{x}}(k+1)=\widetilde{\mathbf{W}}(k)\widetilde{\mathbf{x}}(k)-\alpha_{k}\left(\mathbf{I}-\mathbf{J}\right)\underbrace{\left(\nabla F(\mathbf{x}(k))+\mathbf{v}(k)\right)}_{\text{$\mathbf{w}(k)$}},
\end{align}
where $\widetilde{\mathbf{W}}(k) = \mathbf{W}(k)-\mathbf{J}=\mathbf{I}-\mathbf{L}(k)-\mathbf{J}$. Note that, $\mathbb{E}\left[\left\|\mathbf{w}(k)\right\|^{2}\right]\leq c_{7} < \infty$,
which follows due to the mean square boundedness of $\mathbf{x}(k)$
 and $\nabla F(\mathbf{x}(k))$. Then,  we have:
\begin{align*}
\left\|\widetilde{\mathbf{x}}(k+1)\right\| \le \left\|\widetilde{\mathbf{W}}(k)\right\| \left\|\widetilde{\mathbf{x}}(k)\right\|+\alpha_{k}\left\|{\mathbf{w}}(k)\right\|.
\end{align*}
We now invoke Lemma~4.4 in \cite{soummyaAdaptive} to note that,
%
%
%
after an appropriately chosen $k_{1}$, we have for $\forall k\geq k_1$,
\begin{align}
\label{eq:dis4}
\left\|\widetilde{\mathbf{x}}(k+1)\right\| \le (1-r(k))\left\|\widetilde{\mathbf{x}}(k)\right\|+\alpha_{k}\left\|{\mathbf{w}}(k)\right\|,
\end{align}
with $r(k)$ being a $\mathcal{F}_k$-adapted process that satisfies
$r(k)\in[0,1]$, a.s., and:
\begin{align}
\label{eq:dis3}
\mathbb{E}\left[r(k)|\mathcal{F}_{k}\right]\geq c_8 \beta_{k} = \frac{c_9}{(k+1)^{\frac{1}{2}}}~a.s.,
\end{align}
for some constants $c_8,c_9>0.$
Using \eqref{eq:cool_ineq} in \eqref{eq:dis4}, we have:
\begin{eqnarray*}
&\,& \left\|\widetilde{\mathbf{x}}(k+1)\right\|^{2}\le \left(1+\theta_k\right)(1-r(k))^{2}\left\|\widetilde{\mathbf{x}}(k)\right\|^{2}\\
&+&\left(1+\frac{1}{\theta_k}\right)\alpha_k^2\left\|{\mathbf{w}}(k)\right\|^2,
\end{eqnarray*}
for $\theta_k = \frac{c_{10}}{(k+1)^{\frac{1}{2}}}$.
Then, we have:
\begin{align*}
&\mathbb{E}\left[\left\|\widetilde{\mathbf{x}}(k+1)\right\|^{2}|\mathcal{F}_{k}\right]\leq \left(1+\theta_k\right)\left(1-\frac{c_9}{(k+1)^{\frac{1}{2}}}\right)^{2}\left\|\widetilde{\mathbf{x}}(k)\right\|^{2}\\
&+\left(1+\frac{1}{\theta_k}\right)\alpha_k^2\,
\mathbb{E}[\,\|\mathbf{w}(k)\|^2\,|\,\mathcal{F}_k\,],\,\,\mathrm{a.s.}\nonumber
\end{align*}
%
%
Next, for $c_{10}<c_{9}$~($c_{10}$ can be chosen freely), we have:
\begin{align}
\label{eq:dis7}
&\mathbb{E}\left[\left\|\widetilde{\mathbf{x}}(k+1)\right\|^{2}\right]\leq \left(1-\frac{c_{11}}{(k+1)^{\frac{1}{2}}}\right)
\mathbb{E}\left[\left\|\widetilde{\mathbf{x}}(k)\right\|^{2}\right]\\
&+\frac{c_{12}}{(k+1)^{\frac{3}{2}}}\nonumber
\end{align}
Utilizing Lemma \ref{lemma-estimation},
inequality \eqref{eq:dis7} finally yields
 $\mathbb{E}\left[\left\|\widetilde{\mathbf{x}}(k+1)\right\|^{2}\right] = O\left(\frac{1}{k}\right).$
 The proof of the Lemma is complete.

\subsection{Proof of Theorem~\ref{theorem-1}}
\label{sec:opt_gap}
We are now ready to prove Theorem~\ref{theorem-1}.

\textit{Proof.}

Consider global average $\overline{\mathbf{x}}(k)=\frac{1}{N}\sum_{n=1}\mathbf{x}_{i}(k)$.
From \eqref{eq:dis1}, we have:
\begin{align*}
&\overline{\mathbf{x}}(k+1) = \overline{\mathbf{x}}(k)-\alpha_{k}\left[\frac{1}{N}\sum_{i=1}^{N}\nabla f_{i}\left(\mathbf{x}_{i}(k)\right)+\underbrace{\frac{1}{N}\sum_{i=1}^{N}
\mathbf{v}_{i}(k)}_{\text{$\overline{\mathbf{v}}(k)$}}\right]
\end{align*}
which implies:
\begin{align*}
& \overline{\mathbf{x}}(k+1) = \overline{\mathbf{x}}(k)-\frac{\alpha_k}{N}\left[\sum_{i=1}^{N}\nabla f_{i}\left(\mathbf{x}_{i}(k)\right) \right.\\
&\left. -\nabla f_{i}\left(\overline{\mathbf{x}}(k)\right)+\nabla f_{i}\left(\overline{\mathbf{x}}(k)\right)\right]-\alpha_{k}\overline{\mathbf{v}}(k).
\end{align*}
Recall $f(\cdot)=\sum_{i=1}^{N}f_{i}(\cdot)$.
Then, we have:
\begin{align}
\label{eq:opt2}
&\overline{\mathbf{x}}(k+1) = \overline{\mathbf{x}}(k)-\frac{\alpha_{k}}{N}\nabla f\left(\overline{\mathbf{x}}(k)\right)\\
&-\frac{\alpha_k}{N}\left[\sum_{i=1}^{N}\nabla f_{i}\left(\mathbf{x}_{i}(k)\right)-\nabla f_{i}\left(\overline{\mathbf{x}}(k)\right)\right]-\alpha_{k}\overline{\mathbf{v}}(k),\nonumber
\end{align}
which implies:
\begin{align}
\label{eqn-inexact-centralized}
& \overline{\mathbf{x}}(k+1) = \overline{\mathbf{x}}(k)\\
&-\frac{\alpha_k}{N}\left[\nabla f\left(\overline{\mathbf{x}}(k)\right)+\mathbf{e}(k)\right], \nonumber
\end{align}
where
\begin{align}
\label{eq:opt3}
\mathbf{e}(k) = N\overline{\mathbf{v}}(k)+\underbrace{\sum_{i=1}^{N}\left(\nabla f_{i}\left(\mathbf{x}_{i}(k)\right)-\nabla f_{i}\left(\overline{\mathbf{x}}(k)\right)\right)}_{\text{$\boldsymbol{\epsilon}(k)$}}.
\end{align}
Note that, $\left\|\nabla f_{i}\left(\mathbf{x}_{i}(k)\right)-\nabla f_{i}\left(\overline{\mathbf{x}}(k)\right)\right\| \leq L\left\|\mathbf{x}_{i}(k)-\overline{\mathbf{x}}(k)\right\| = L\left\|\widetilde{\mathbf{x}}_{i}(k)\right\|$.
Thus, we can conclude for
\begin{align*}
&\boldsymbol{\epsilon}(k) = \sum_{i=1}^{N}\left(\nabla f_{i}\left(\mathbf{x}_{i}(k)\right)-\nabla f_{i}\left(\overline{\mathbf{x}}(k)\right)\right)
\end{align*}
the following:
\begin{align}
\label{eqn-bias-grad}
%
 \mathbb{E}\left[\left\|\boldsymbol{\epsilon}(k)\right\|^{2}\right] \leq \frac{c_{15}}{(k+1)}.
\end{align}
Note here that \eqref{eqn-inexact-centralized}
is an inexact gradient method for minimizing~$f$ with step size
$\alpha_k/N$
and the random gradient error
$\mathbf{e}(k)=N \overline{\mathbf{v}}(k)+\boldsymbol{\epsilon}(k)$.
The term $N \overline{\mathbf{v}}(k)$ is zero-mean, while
the gradient estimate bias is induced by $\boldsymbol{\epsilon}(k)$;
as per~\eqref{eqn-bias-grad}, the bias is at most $O(1/k)$ in the mean square sense.

With the above development in place, we rewrite \eqref{eq:opt2} as follows:
\begin{align}
\label{eq:opt5}
&\overline{\mathbf{x}}(k+1) = \overline{\mathbf{x}}(k)-\frac{\alpha_{k}}{N}\nabla f\left(\overline{\mathbf{x}}(k)\right)-\frac{\alpha_k}{N}\boldsymbol{\epsilon}(k)
-\alpha_{k}\overline{\mathbf{v}}(k).
\end{align}
This implies, recalling that $\mathbf{x}^\star$ is the solution to \eqref{eq:opt_problem}:
\begin{align}
& \overline{\mathbf{x}}(k+1)-{\mathbf{x}^{\star}} = \overline{\mathbf{x}}(k)-{\mathbf{x}^{\star}}\\
& -\frac{\alpha_{k}}{N}\left[\nabla f\left(\overline{\mathbf{x}}(k)\right)-\underbrace{\nabla f\left({\mathbf{x}^{\star}}\right)}_{\text{$=0$}}\right]-\frac{\alpha_k}{N}\boldsymbol{\epsilon}(k)-\alpha_{k}\overline{\mathbf{v}}(k).
\end{align}
By the mean value theorem, we have:
\begin{align}
& \nabla f\left(\overline{\mathbf{x}}(k)\right)-\nabla f\left({\mathbf{x}^{\star}}\right) = \underbrace{\left[\int_{s=0}^{1}\nabla^{2}
f\left({\mathbf{x}^{\star}}+s\left(\overline{\mathbf{x}}(k)
-{\mathbf{x}^{\star}}\right)\right)\right]\,d\,s}_{\text{$\overline{\mathbf{H}}_{k}$}}
 \nonumber \\
 \label{eq:opt6}
& \times \left(\overline{\mathbf{x}}(k)-{\mathbf{x}^{\star}}\right),
\end{align}
where it is to be noted that $NL\succcurlyeq\overline{\mathbf{H}}_{k}\succcurlyeq N\mu$.
Using \eqref{eq:opt6} in  \eqref{eq:opt5}, we have:
\begin{align}
\label{eq:opt7}
&\left(\overline{\mathbf{x}}(k+1)-{\mathbf{x}^{\star}}\right) = \left[\mathbf{I}-\frac{\alpha_k}{N}\overline{\mathbf{H}}_{k}\right]\left(\overline{\mathbf{x}}(k) -{\mathbf{x}^{\star}}\right)\\
&-\frac{\alpha_k}{N}\boldsymbol{\epsilon}(k)-\alpha_{k}\overline{\mathbf{v}}(k). \nonumber
\end{align}
Denote by $\mathbf{m}(k)=\left[\mathbf{I}-
\frac{\alpha_k}{N}\overline{\mathbf{H}}_{k}\right]\left(\overline{\mathbf{x}}(k)
-{\mathbf{x}^{\star}}\right)-\frac{\alpha_k}{N}\boldsymbol{\epsilon}(k)$.
Then, \eqref{eq:opt7} is rewritten as:
\begin{align}
\label{eq:opt8}
&\left(\overline{\mathbf{x}}(k+1)-{\mathbf{x}^{\star}}\right)=\mathbf{m}(k)-\alpha_{k}\overline{\mathbf{v}}(k),
\end{align}
and so:
\begin{align*}
&  \left\|\overline{\mathbf{x}}(k+1)-{\mathbf{x}^{\star}}\right\|^{2} \le \left\|\mathbf{m}(k)\right\|^{2}-2\alpha_{k}\mathbf{m}(k)^{\top}\overline{\mathbf{v}}(k)\\
&+\alpha_{k}^{2}\left\|\overline{\mathbf{v}}(k)\right\|^{2}.
\end{align*}
The latter inequality implies:
\begin{align*}
&  \mathbb{E}[\,\left\|\overline{\mathbf{x}}(k+1)-{\mathbf{x}^{\star}}\right\|^{2}\,|\,\mathcal{F}_k\,] \le \left\|\mathbf{m}(k)\right\|^{2}\\
&-2\alpha_{k}\mathbf{m}(k)^{\top}\mathbb{E}[\,\overline{\mathbf{v}}(k)\,|\,\mathcal{F}_k\,]+
\alpha_{k}^{2}\mathbb{E}[\,\left\|\overline{\mathbf{v}}(k)\right\|^{2}\,|\,\mathcal{F}_k],\,\,\mathrm{a.s.}
\end{align*}
Taking expectation, using the fact that $\mathbb{E}[\,\overline{\mathbf{v}}(k)\,|\,\mathcal{F}_k\,]=0$,
Assumption~\ref{assumption-gradient-noise}, and Lemma~\ref{Lemma-MSS-BDD},
we obtain:
\begin{align}
\label{eq:opt9}
\mathbb{E}\left[\left\|\overline{\mathbf{x}}(k+1)-{\mathbf{x}^{\star}}\right\|^{2}\right] \le \mathbb{E}\left[\left\|\mathbf{m}(k)\right\|^{2}\right]+\frac{c_{17}}{(k+1)^2},
\end{align}
for some constant $c_{17}>0.$
Next, using \eqref{eq:cool_ineq}, we have for $\mathbf{m}(k)$ the following:
\begin{align}
&\left\|\mathbf{m}(k)\right\|^{2} \le \left(1+\theta_{k}\right)\left\|\mathbf{I}-
\frac{\alpha_k}{N}\overline{\mathbf{H}}_{k}\right\|^{2}\left\|\overline{\mathbf{x}}(k)-{\mathbf{x}^{\star}}\right\|^{2}
\nonumber \\
&+\left(1+\frac{1}{\theta_k}\right)\frac{\alpha_k^2}{N^2}\left\|\boldsymbol{\epsilon}(k)\right\|^{2}\nonumber\\
&\le \left(1+\theta_{k}\right)(1-c_{18}\alpha_{k})^{2}\left\|\overline{\mathbf{x}}(k)-{\mathbf{x}^{\star}}\right\|^{2}
\nonumber \\
&+\left(1+\frac{1}{\theta_k}\right)\frac{\alpha_k^2}{N^2}\left\|\boldsymbol{\epsilon}(k)\right\|^{2}, \nonumber
\end{align}
with $c_{18}=\mu/N$, because $\mu\,\mathbf{I} \preceq \overline{\mathbf{H}}_{k} \preceq L\,\mathbf{I}$.
 After choosing $\theta_{k}=\frac{c_{19}}{(k+1)}$ such that $c_{19} < \alpha_0\,c_{18}/2 =
 \alpha_0\,\mu/(2N)$
 and after taking expectation, we obtain:
\begin{align}
\label{eq:opt11}
&\mathbb{E}[\,\left\|\mathbf{m}(k)\right\|^{2} \,]\leq \left(1-\frac{c_{20}}{k+1}\right)\mathbb{E}[ \, \left\|\overline{\mathbf{x}}(k)
-{\mathbf{x}^{\star}}\right\|^{2}\,]+\frac{c_{21}}{(k+1)^2},
\end{align}
where
$c_{20}>\alpha_0\,\mu/(2N)>1$ (because
$\alpha_0 > 2N/\mu$)  and $c_{21}$ is a positive constant.
Combining \eqref{eq:opt11} and
\eqref{eq:opt9}, we get:
\begin{align*}
& \mathbb{E}\left[\left\|\overline{\mathbf{x}}(k+1)-{\mathbf{x}^{\star}}\right\|^{2}\right] \leq \left(1-\frac{c_{20}}{k+1}\right)\left\|\overline{\mathbf{x}}(k)-{\mathbf{x}^{\star}}\right\|^{2}\\
& +\frac{c_{21}}{(k+1)^2}+\frac{c_{17}}{(k+1)^2}.
\end{align*}
Invoking Lemma~\ref{lemma-estimation}, the latter inequality implies:
\begin{align}
\label{eqn-invoking}
\mathbb{E}\left[\left\|\overline{\mathbf{x}}(k+1)-{\mathbf{x}^{\star}}\right\|^{2}\right] \leq \frac{c_{22}}{(k+1)},
\end{align}
for some constant $c_{22}>0.$
Therefore, for the global average $\overline{\mathbf{x}}(k)$, we have obtained the mean square rate
$O\left(\frac{1}{k}\right)$.
Finally, we note that,
\begin{align}
\label{eq:opt12}
&\left\|\mathbf{x}_{i}(k)-{\mathbf{x}^{\star}}\right\| \le \left\|\overline{\mathbf{x}}(k)-{\mathbf{x}^{\star}}\right\|+
\left\|\underbrace{\mathbf{x}_{i}(k)-\overline{\mathbf{x}}(k)}_{\text{$\widetilde{\mathbf{x}}_{i}(k)$}}\right\|.\nonumber
\end{align}
After using:
\begin{align}
&\left\|\mathbf{x}_{i}(k)-{\mathbf{x}^{\star}}\right\|^{2} \leq 2\left\|\widetilde{\mathbf{x}}_{i}(k)\right\|^{2}+
2\left\|\overline{\mathbf{x}}(k)-{\mathbf{x}^{\star}}\right\|^{2},\nonumber
\end{align}
and taking expectation,
it follows that $\mathbb{E}\left[\left\|\mathbf{x}_{i}(k)-{\mathbf{x}^{\star}}\right\|^{2}\right] = O\left(\frac{1}{k}\right)$, for all $i=1,...,N$. The proof is complete.
%

\section{Simulation example}

We provide a simulation example
on $\ell_2$-regularized logistic losses and
random networks where links fail independently
over iterations and across different links, with probability $p_{\mathrm{fail}}$.
The simulation corroborates the derived $O(1/k)$
rate of algorithm \eqref{eqn-alg-node-i} over random networks and shows that
deterioration due to increase of $p_{\mathrm{fail}}$ is small.

We consider empirical risk minimization \eqref{eqn-empir-risk-fcn}
with the logistic loss in \eqref{eqn-ML-examples} and
the regularization functions set to $\Psi_i(\mathbf{x})=
\frac{\kappa}{2}\|\mathbf{x}\|^2$, $i=1,...,N$,
where $\kappa>0$ is the regularization parameter that
is set to $\kappa=0.5$.

The number of data points per node is
$n_i=10$. We generate the ``true'' classification vector $x^\prime=((\mathbf{x}_1^\prime)^\top, x_0^\prime)^\top$
by drawing its entries independently from standard normal distribution.
Then, the class labels are generated as $b_{ij}=\mathrm{sign} \left( (\mathbf{x}^\prime_1)^\top \mathbf{a}_{i,j}+x^\prime_0+\epsilon_{ij}\right)$, where
$\epsilon_{ij}$'s are drawn independently from normal distribution with zero mean and
standard deviation~$2$.
The feature vectors $\mathbf{a}_{i,j}$, $j=1,...,n_i$, at node $i$
are generated as follows: each entry of
each vector is a sum of
a standard normal random variable
and a uniform random variable with support
$[0,\,5\,i]$. Different entries within a feature vector
are drawn independently, and also different vectors are drawn independently,
both intra node and inter nodes. Note that
the feature vectors at different nodes
are drawn from different distributions.

The algorithm parameters are set as follows.
 We let $\beta_k=\frac{1}{\theta\,(k+1)^{1/2}}$,
 $\alpha_k=\frac{1}{k+1}$, $k=0,1,...$
 Here, $\theta$ is the maximal
 degree across all nodes in the network
 and here equals $\theta=6$.
  Algorithm \eqref{eqn-alg-node-i}
  is initialized with
  $\mathbf{x}_i(0)=0$, for all
  $i=1,...,N$.

We consider a connected network $\mathcal{G}$ with $N=10$ nodes and $23$ links,
 generated as a random geometric graph: nodes are placed
 randomly (uniformly) on a unit square, and the node pairs whose distance is less
than a radius are connected by an edge.
 We consider
 the random network model
 where each (undirected) link
 in network $\mathcal{G}$ fails
 independently across iterations and
 independently from other links
 with probability~$p_{\mathrm{fail}}$.
 We consider the cases $p_{\mathrm{fail}} \in \{0;\,0.5;\,0.9\}$.
 Note that the case $p_{\mathrm{fail}}=0$
 corresponds to network~$\mathcal{G}$
 with all its links always online, more precisely,
 with links failing with zero probability.
 Algorithm~\eqref{eqn-alg-node-i} is
  then run on each of the described network models, i.e.,
  for each $p_{\mathrm{fail}} \in \{0;\,0.5;\,0.9\}$.
  This allows us to assess how much
  the algorithm performance degrades with the increase
  of~$p_{\mathrm{fail}}$.
    We also include a
 comparison with the following centralized stochastic
 gradient method:
 \begin{equation}
 \label{eqn-centralized-SGD}
 \mathbf{y}(k+1) = \mathbf{y}(k) -
 \frac{1}{N(k+1)}
 \sum_{i=1}^N \nabla \ell\left(\,\mathbf{y}(k);\,\mathbf{a}_i(k),b_i(k)\,\right),
 \end{equation}
 where $(\mathbf{a}_i(k),b_i(k))$
  is drawn uniformly
  from the set $(\mathbf{a}_{i,j},b_{i,j})$, $j=1,...,n_i$.
  Note that algorithm \eqref{eqn-centralized-SGD}
  makes an unbiased estimate
  of $\sum_{i=1}^N \nabla f_i(\mathbf{y}(k))$
  by drawing a sample uniformly at random
  from each node's data set.
  Algorithm~\eqref{eqn-centralized-SGD} is an idealization of~\eqref{eqn-alg-node-i}:
   it shows how~\eqref{eqn-alg-node-i} would be implemented
   if there existed a fusion node that had access to all nodes' data.
   Hence, the comparison with~\eqref{eqn-centralized-SGD}
   allows us to examine
    how much the performance of \eqref{eqn-alg-node-i} degrades due to lack
    of global information, i.e., due to the distributed nature of the considered problem.
   Note that step-size in~\eqref{eqn-centralized-SGD}
   is set to $1/N(k+1)$ for a meaningful comparison with~\eqref{eqn-alg-node-i},
   as this is the step-size effectively utilized by
   the hypothetical global average of the nodes' iterates with~\eqref{eqn-alg-node-i}.
As an error metric, we use
   the mean square error (MSE) estimate averaged across nodes:
    $
   \frac{1}{N} \sum_{i=1}^N {\|\mathbf{x}_i(k)-{\mathbf{x}}^\star\|^2}.
   $

Figure~1 plots the estimated MSE, averaged
across 100 algorithm runs,
versus iteration number~$k$ for
different values of parameter~$p_{\mathrm{fail}}$
in $\mathrm{log}_{10}$-$\mathrm{log}_{10}$ scale.
Note that here the slope of the plot curve
corresponds to the sublinear rate of the method;
e.g., the $-1$ slope corresponds to a $1/k$ rate.
 First, note from the Figure
 that, for any value of~$p_{\mathrm{fail}}$,
 algorithm \eqref{eqn-alg-node-i}
 achieves on this example (at least) the $1/k$
 rate, thus corroborating our theory.
 Next, note that
 the increase of the link failure probability
 only increases the constant in the MSE but does not
 affect the rate. (The curves
 that correspond to different values
 of $p_{\mathrm{fail}}$ are
 ``vertically shifted.'')
  Interestingly, the loss due to
  the increase of $p_{\mathrm{fail}}$ is
  small; e.g., the curves
  that correspond to $p_{\mathrm{fail}}=0.5$ and
  $p_{\mathrm{fail}}=0$ (no link failures) practically match.
Figure~1 also shows the performance of
the centralized method~\eqref{eqn-centralized-SGD}.
We can see that, except for
the initial few iterations, the distributed method \eqref{eqn-alg-node-i}
is very close in performance to the centralized method.

\begin{figure}[thpb]
      \centering
      \includegraphics[height=3.4 in,width=2.8 in, angle =-90]{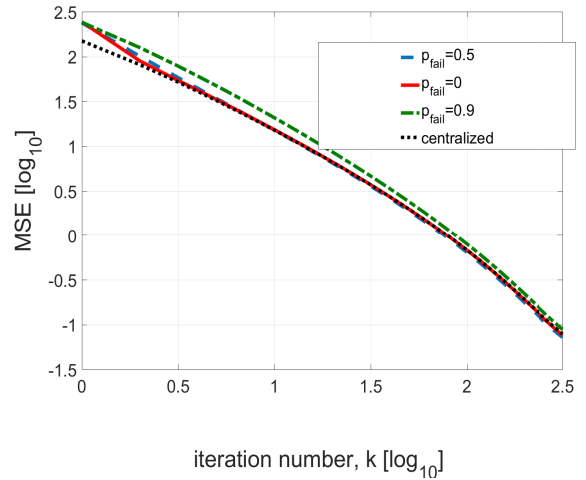}
      \caption{Estimated MSE versus iteration number~$k$ for algorithm \eqref{eqn-alg-node-i}
      with link failure probability $p_{\mathrm{fail}}=0$ (red, solid line);
      $0.5$ (blue, dashed line); and $0.9$ (green, dash-dot line).
      The Figure also shows the performance of
      the centralized stochastic gradient method in~\eqref{eqn-centralized-SGD}
      (black, dotted line).}
      \label{Figure_1}
      \vspace{-9mm}
\end{figure}

\vspace{3mm}
\section{Conclusion}
We considered a distributed stochastic gradient method
for smooth strongly convex optimization.
Through the analysis of the considered method,
we established for the first time the order optimal
$O(1/k)$ MSE convergence rate for the assumed
optimization setting when the underlying network is
randomly varying. Simulation example
on $\ell_2$-regularized logistic losses
corroborates the established $O(1/k)$
rate and suggests that the effect of the
underlying network's statistics on
the $O(1/k)$ rate's hidden constant is small.

\bibliographystyle{IEEEtran}
\bibliography{IEEEabrv,bibliographyCDC2018}

\begin{thebibliography}{10}
\providecommand{\url}[1]{#1}
\csname url@samestyle\endcsname
\providecommand{\newblock}{\relax}
\providecommand{\bibinfo}[2]{#2}
\providecommand{\BIBentrySTDinterwordspacing}{\spaceskip=0pt\relax}
\providecommand{\BIBentryALTinterwordstretchfactor}{4}
\providecommand{\BIBentryALTinterwordspacing}{\spaceskip=\fontdimen2\font plus
\BIBentryALTinterwordstretchfactor\fontdimen3\font minus
  \fontdimen4\font\relax}
\providecommand{\BIBforeignlanguage}[2]{{%
\expandafter\ifx\csname l@#1\endcsname\relax
\typeout{** WARNING: IEEEtran.bst: No hyphenation pattern has been}%
\typeout{** loaded for the language `#1'. Using the pattern for}%
\typeout{** the default language instead.}%
\else
\language=\csname l@#1\endcsname
\fi
#2}}
\providecommand{\BIBdecl}{\relax}
\BIBdecl

\bibitem{SoummyaGossipEstimation}
S.~Kar and J.~M.~F. Moura, ``Convergence rate analysis of distributed gossip
  (linear parameter) estimation: {F}undamental limits and tradeoffs,''
  \emph{IEEE Journal of Selected Topics in Signal Processing, Signal Processing
  in Gossiping Algorithms Design and Applications}, vol.~5, no.~4, pp.
  674–--690, Aug. 2011.

\bibitem{BulloBook}
F.~Bullo, J.~Cortes, and S.~Martinez, \emph{Distributed control of robotic
  networks: {A} mathematical approach to motion coordination algorithms}.\hskip
  1em plus 0.5em minus 0.4em\relax Princeton University Press, 209.

\bibitem{daneshmand2015hybrid}
A.~Daneshmand, F.~Facchinei, V.~Kungurtsev, and G.~Scutari, ``Hybrid
  random/deterministic parallel algorithms for convex and nonconvex big data
  optimization,'' \emph{IEEE Transactions on Signal Processing}, vol.~63,
  no.~15, pp. 3914--3929, 2015.

\bibitem{AsuRandom2}
I.~Lobel and A.~E. Ozdaglar, ``Distributed subgradient methods for convex
  optimization over random networks,'' \emph{IEEE Trans. Automat. Contr.},
  vol.~56, no.~6, pp. 1291--1306, Jan. 2011.

\bibitem{ASU_Math_Prog}
I.~Lobel, A.~Ozdaglar, and D.~Feijer, ``Distributed multi-agent optimization
  with state-dependent communication,'' \emph{Mathematical Programming}, vol.
  129, no.~2, pp. 255--284, 2011.

\bibitem{randomNesterov}
D.~Jakovetic, J.~Xavier, and J.~M.~F. Moura, ``Convergence rates of distributed
  {N}esterov-like gradient methods on random networks,'' \emph{IEEE
  Transactions on Signal Processing}, vol.~62, no.~4, pp. 868--882, February
  2014.

\bibitem{RabbatDistributedStronglyCVX}
K.~Tsianos and M.~Rabbat, ``Distributed strongly convex optimization,''
  \emph{50th Annual Allerton Conference onCommunication, Control, and
  Computing}, Oct. 2012.

\bibitem{SayedStochasticOpt}
Z.~J. Towfic, J.~Chen, and A.~H. Sayed, ``Excess-risk of distributed stochastic
  learners,'' \emph{IEEE Transactions on Information Theory}, vol.~62, no.~10,
  Oct. 2016.

\bibitem{DistributedMirrorDescent}
D.~Yuan, Y.~Hong, D.~W.~C. Ho, and G.~Jiang, ``Optimal distributed stochastic
  mirror descent for strongly convex optimization,'' \emph{Automatica},
  vol.~90, pp. 196--203, April 2018.

\bibitem{Kozat}
N.~D. Vanli, M.~O. Sayin, and S.~S. Kozat, ``Stochastic subgradient algorithms
  for strongly convex optimization over distributed networks,'' \emph{IEEE
  Transactions on network science and engineering}, vol.~4, no.~4, pp.
  248--260, Oct.-Dec. 2017.

\bibitem{NedicStochasticPush}
A.~Nedic and A.~Olshevsky, ``Stochastic gradient-push for strongly convex
  functions on time-varying directed graphs,'' \emph{IEEE Transactions on
  Automatic Control}, vol.~61, no.~12, pp. 3936--3947, Dec. 2016.

\bibitem{CDCKW2018}
A.~K. Sahu, D.~Jakovetic, D.~Bajovic, and S.~Kar, ``Distributed zeroth order
  optimization over random networks: {A} {K}iefer-{W}olfowitz stochastic
  approximation approach,'' 2018, available at
  \url{https://www.dropbox.com/s/kfc2hgbfcx5yhr8/MainCDC2018KWSA.pdf}.

\bibitem{soummyaAdaptive}
S.~Kar, J.~M.~F. Moura, and H.~V. Poor, ``Distributed linear parameter
  estimation: {A}symptotically efficient adaptive strategies,'' \emph{SIAM J.
  Control and Optimization}, vol.~51, no.~3, pp. 2200--2229, 2013.

\end{thebibliography}

%
%
%
%

\end{document}